\documentclass[lettersize,journal]{IEEEtran}
\usepackage{amsmath,amsfonts,amssymb,amsthm}
\usepackage{algorithmic}
\usepackage{algorithm}
\usepackage{array}
\usepackage[caption=false,font=normalsize,labelfont=sf,textfont=sf]{subfig}
\usepackage{textcomp}
\usepackage{stfloats}
\usepackage{url}
\usepackage{verbatim}
\usepackage{graphicx}
\usepackage{cite}

\usepackage{nicefrac}
\usepackage{xcolor}
\usepackage{tabularx}
\usepackage{booktabs}
\usepackage{ragged2e}
\usepackage{siunitx}
\usepackage{multirow}
\usepackage{hyperref}

\newtheorem{theorem}{Theorem}

\newtheorem{example}[theorem]{Example}
\newcommand{\R}{\mathbb{R}}

\newcommand{\C}{\mathbb{C}}
\newcommand{\eps}{\varepsilon}

\newcommand{\dlatent}{d_{\mathrm{latent}}}
\newcommand{\dhidden}{d_{\mathrm{hidden}}}

\begin{document}

\title{An Encoder-Transformer Architecture for Recognition of the
Jordan  Structure  of a~Matrix}

\author{Michał Trojanowski, Michał Wojtylak}

\markboth{}
{}
\IEEEpubid{ }

\maketitle

\begin{abstract}

We propose a~machine-learning framework for detecting whether a~given matrix is a~perturbation of a~matrix with a~large Jordan block. The proposed model achieves high classification accuracy on synthetically generated, robustly perturbed data and outperforms a~classical numerical baseline. Moreover, we demonstrate that the learned model generalizes to several classes of matrices not seen during training. These results suggest that the architecture captures structural properties associated with matrix defectiveness.
\end{abstract}


\section{Introduction}

\IEEEPARstart{C}{}omputing the Jordan form of a~nondiagonalisable matrix in finite precision has been always a~challenge for numerical algorithms, see Section~\ref{sec:MB}  for a~short overview. The round-off error makes the matrix always diagonalizable, obscuring the true structure. This causes several problems in mathematical modelling, as sometimes the dynamics follow  the pattern of the original matrix rather than its backward error update \cite{Wojtylak2026}.
In this work we propose a~novelty approach. Driven by the algebraic properties of nilpotency, our architecture processes a~sequence of successive matrix powers using a~combination of dimension-specific multilayer perceptron encoders and transformer network.

\subsection{Problem formulation}

A question that we consider can  be intuitively formulated as:

\emph{ What is the largest size $m$ of a~Jordan block near the matrix~$A$?}

A  precise mathematical formulation requires the following notation: $\| \cdot \|$ denotes the operator 2-norm on $\C^{d\times d}$ and by $\kappa(A):=\|A\|\|A^{-1}\|$ the condition number of a~nonsingular matrix $A$. 
The spectral radius is denoted by $\rho(A)$.
Our question can can be reformulated now as follows.

\emph{Given a~matrix $A \in \C^{d\times d}$  with spectral radius $\rho(A) \leqslant \rho_{\max}$, and $\eps_{\max}>0$, determine the largest integer $m$ such that 
there exists a~Jordan matrix $J$ with zero eigenvalue only, with the largest block of size $m$, and an invertible matrix  $S$ satisfying 
\begin{equation}\label{eq:constraints}
\kappa(S) \leqslant \kappa_{\max} ,\quad \text{ and }\left\|S^{-1}AS - J\right\| \leqslant \eps_{\max}.
\end{equation}
}

Later on we will say in such situation that \emph{the Jordan block size $m$ is \textit{attainable} under perturbation of magnitude $\eps_{\max}$}. 
For practical use it is recommended to use $\eps_{\max}=10^{-1}\rho_{\max}$. 
The value of $\rho_{\max}$ is a~size of a~typical cluster of eigenvalues, typically $10^{-6}$. For numerical stability purposes the current model was trained and tested with $\rho_{\max}=1$.  The value of the constraint $\kappa_{\max}$ is set to $ 200 d$, which  is motivated by asymptotic conditioning of normally distributed matrices, see Theorem \ref{theorem:Edelman}.
A combination of the two constraints in \eqref{eq:constraints} represents the formalization of the heuristic question posed at the beginning of this section, cf. Subsection~\ref{section:discussion} for a~further discussion.

This approach also allows to compute the Jordan form of a~matrix with several clusters of eigenvalues, using the Riesz projection approach, see Section~\ref{section:nilpotent} for details.

\subsection{Motivation I: Dynamical systems}

Consider an autonomous system of ODEs $x'(t) = Ax(t)$, $x(0)=x_0$ with $A\in\C^{d\times d}$. While the eigenvalues of $A$~determine the asymptotic growth or decay of the solution, the Jordan canonical form of $A$ dictates structural stability. The presence of a~Jordan block of size $k$ associated with eigenvalue $\lambda$ introduces terms of the form $t^je^{\lambda t}$, $j=0,...,k-1$ in the solution $x(t)$. Even if $\mathrm{Re}(\lambda)<0$ ensures eventual stability, a~large value of $k$ can cause significant growth before the exponential decay. When $\mathrm{Re}(\lambda)=0$, the defectivity of $\lambda$~becomes the deciding factor for stability.

If a~matrix $A$ lies in the neighbourhood of a~matrix with a~large Jordan block, it implies that the system's behaviour is highly sensitive to noise and susceptible to instability. Therefore, for a~matrix with a~cluster of eigenvalues, we are interested in the maximal attainable Jordan block size.

\subsection{Motivation II: Image recognition}
Neural networks, including CNNs and Vision Transformers, excel at tasks like image classification\cite{Krizhevsky2012, Kolesnikov2021} or denoising \cite{Vincent2010, Kulkarni2023}. Image recognition architectures rely heavily on the spatial local correlation of the tensor representing the image, such as convolutional structures in CNNs \cite{Krizhevsky2012} or a~sequence of patches in ViTs \cite{Kolesnikov2021}. 

\begin{figure}[tbp]
    \centering
    \includegraphics[width=\linewidth]{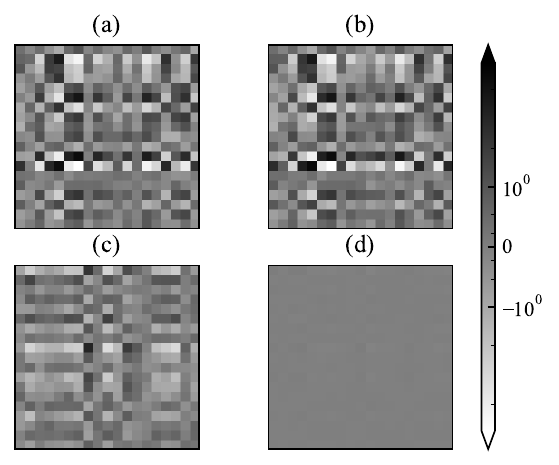}
    \caption{Visualization (asinh grayscale) of four matrices in the neighbourhood of a~matrix with a~single eigenvalue $\lambda=0$ of algebraic multiplicity 19. Despite their visual similarity, matrices (a) and (b) correspond to different structural properties: (a) is a~perturbation of a~matrix with a~maximal Jordan block of size 4, whereas (b) corresponds to a~maximal block of size 8. Conversely, matrix (c) appears visually distinct despite also having a~maximal block size of 4. Finally, (d) represents a~perturbed zero matrix (diagonalizable), which can be easily recognized. Our model correctly predicted all maximal attainable block sizes.}

    \label{fig:random_matrices}
\end{figure}

However, the Jordan structure of a~matrix presents a~distinct challenge: it is invariant under similarity transformations. Consequently, matrix entries lack the local dependence, making standard vision models poorly suited for direct application. The difficulty of identifying the maximal Jordan block size is exemplified in Figure \ref{fig:random_matrices}: patterns are often invisible to the naked eye or even illusory, making automated recognition both necessary and challenging. We will compare the difficulty of learning Jordan structure recognition by the model to image recognition tasks.  

Despite the lack of spatial structure, we draw inspiration from the ViT flow \cite{Kolesnikov2021}. Instead of extracting features based on visual-semantic regions, our model utilizes algebraic features derived from linear algebra theory. This shift allows the attention mechanism to focus on the global structural properties of the operator rather than local entry-wise correlations.

\subsection{Existing research}

Demmel in \cite{Demmel1983}  initiated the study of finding the nearest defective matrix to a~given diagonalizable matrix. His work was continued by Alam et al. in \cite{Alam2005, Alam2011}. They proved existence of a~(not necessarily unique) nearest defective matrix and proposed a~numerical construction based on Newton iteration. 

While this problem initially appears closely linked to ours, the relationship is subtle. Consider the example of a~perturbed Jordan block of size 3 from \cite[Section 6]{Alam2011}. It appears that despite the proximity of a~larger Jordan block, the nearest defective matrix has two distinct eigenvalues: one simple and one with defect 2. The authors conjecture that under generic assumptions, the algebraic multiplicity of the defective eigenvalue in the nearest defective matrix is exactly 2.~Consequently, in the generic case, it is impossible to obtain Jordan blocks larger than size 2 (or eigenvalues with a~defect greater than 1) by computing the nearest defective matrix.

\section{Mathematical Analysis}\label{sec:MB}

A matrix $A \in \C^{d \times d}$ is called \textit{diagonalizable} if  it can be represented as $A = SDS^{-1}$ for some diagonal matrix $D$~and invertible $S$. Otherwise, the matrix is \textit{defective}. 
We say an eigenvalue $\lambda$ of $A$ is \textit{multiple} if its algebraic multiplicity is greater than one, and \textit{defective} if its algebraic multiplicity strictly exceeds its geometric multiplicity.

\subsection{Problem setting discussion}\label{section:discussion}

Let us explain now in detail the mathematical formulation of the problem.
For a~matrix $A \in \C^{d\times d}$, let $\delta_k(A)$ denote the distance to the nearest matrix whose Jordan canonical form has a Jordan block of size $k$. 
This distance changes under similarity transformations by at most the condition number of the similarity matrix. More precisely, for every nonsingular $S \in \C^{d\times d}$,
\begin{equation}\label{eq:lemma}
    \frac{\delta_k(A)}{\kappa(S)}
    \leqslant
    \delta_k(SAS^{-1})
    \leqslant
    \kappa(S)\,\delta_k(A),
\end{equation}
  see \cite[Lemma 2.1]{Demmel1983}.
In particular,
$
\delta_k(A)
    \leqslant
    {\kappa_{\max}}\,\eps_{\max}.
$

Further, notice that by a simple Schur-form argument 
$\delta_k(A)$ is estimated from above by the diameter of the smallest cluster of $k$ eigenvalues. 
Hence, computing $\delta_k(A)$ 
would trivialize in practice to the problem to determining the diameter of the  cluster. This is not the aim, as such nearby defective matrix would usually have extreme values of the conditioning of the Jordan basis. Therefore, two parallel conditions in \eqref{eq:constraints} are imposed.

\begin{example}
    Consider the matrix 
\begin{equation*}
    A_\delta = \begin{pmatrix}
    10^{-3} & 1 & 0 & 0 \\
    0 & 10^{-2} & \delta & 0 \\
    0 & 0 & 0 & 0 \\
    0 & 0 & 0 & 0
    \end{pmatrix}
\end{equation*}
with a~cluster of eigenvalues at the origin. While $A_\delta$ may seem close to a~matrix with Jordan block of size 3 for any $\delta \neq 0$, the Jordan basis of the neighbouring matrix has condition number of order $O(\delta^{-1})$. Consequently, under the constraint of a~bounded condition number, the matrix is effectively limited to a~maximal attainable block size of 2, as seen in its $2 \times 2$ upper-left structure.
\end{example}

\subsection{Matrix Conditioning and Eigenvalues}
 We recall several results relating matrix conditioning, eigenvalue stability, and the typical behavior of random matrices.
All these results justify the choice of the values of the parameters $\rho_{\max}$, $\eps_{\max}$, and $\kappa_{\max}$.

For $A$ diagonalizable, let 
\begin{equation}
    \kappa_V(A) := \inf\{\kappa(V) : \, A = V^{-1}DV,\; D \text{ diagonal}\}
\end{equation}
denote the eigenvector condition number of $A$. The Bauer--Fike theorem \cite{Bauer1960} is an important result showing the stability of matrix eigenvalues under good eigenvector conditioning.

\begin{theorem}
    For $A\in\C^{d\times d}$ diagonalizable, $E\in\C^{d\times d}$, and $\tilde{\lambda}$ eigenvalue of $A + E$, there exists eigenvalue $\lambda$ of $A$ such that
    \begin{equation}
        \left|\tilde{\lambda} - \lambda\right| \leqslant \kappa_V(A)\|E\|.
    \end{equation}
\end{theorem}

We now cite two theorems regarding the regularization of eigenvalue conditioning. The first, by Banks et al. \cite[Theorem 1.1]{Banks2021}, states that any complex matrix can be well approximated by a~matrix with low eigenvector conditioning; thus, poor eigenvector conditioning is not structurally stable and can be improved by arbitrarily small perturbations.
\begin{theorem}
    For $A\in\C^{d\times d}$ and $\delta\in(0,1)$, there exists a~matrix $E\in\C^{d \times d}$ such that $\|E\|\leqslant\delta\|A\|$ and
    \begin{equation}
        \kappa_V(A+E) \leqslant 4d^{\frac{3}{2}}\left(1+\delta^{-1}\right).
    \end{equation}
\end{theorem}

The second result, by Blazhko and Wojtylak \cite{Wojtylak2026}, provides a~probabilistic perspective on this phenomenon.

\begin{theorem}
    Let $A\in \C^{d\times d}$ and let $E$ be a~complex Ginibre matrix, i.e., has i.i.d. complex entries with both real and imaginary parts sampled from $\mathcal{N}\left(0, 1/{2d}\right)$. Let $\eps \in \left(0, \|A\|\right)$. Then:
    \begin{equation}
         \mathbb{E}\kappa_V\left(A + \eps E\right) \leqslant \frac{\alpha(d)\|A\|}{\eps},
    \end{equation}
    where $\alpha(d)$ is a~constant depending only on $d$.
\end{theorem}

At the end of the theoretical introduction, the following result, due to Edelman \cite[Theorem 7.1]{Edelman1989}, quantifies the asymptotic growth rate and distribution of the condition number of a~random matrix.
\begin{theorem}\label{theorem:Edelman}
    Let $A \in \R^{d\times d}$ be a~random matrix with i.i.d. entries from $\mathcal{N}(0,1)$. Then the PDF of $\kappa(A)/d$ converges pointwise to
    \begin{equation}
        f(x) = \frac{2x+4}{x^3}e^{-\frac{2}{x}-\frac{2}{x^2}}
    \end{equation}
    as $d \rightarrow \infty$. Moreover, $\mathbb{E}\left[\log\kappa(A)\right] = \log d + c + o(1)$ with $c \approx 1.537$ as $d\rightarrow\infty$.
\end{theorem}

\subsection{Reduction to nilpotent matrices}\label{section:nilpotent}

For a~matrix $A \in \C^{n\times n}$, we  denote its spectrum by $\Lambda(A)$. Consider a~specific cluster of eigenvalues $\Lambda_0 \subset \Lambda(A)$ with total algebraic multiplicity $d$. Let $\Gamma$ be a~simple, positively oriented, closed contour enclosing a~region $G_\Gamma$ such that $G_\Gamma \cap \Lambda(A) = \Lambda_0$. We assume that $\Gamma \cap \Lambda(A) = \emptyset$.  Under these conditions, we define the Riesz projection:

\begin{equation}\label{eq:Riesz_projection}
    P_\Gamma := \frac{1}{2\pi i} \oint_\Gamma \left(zI - A\right)^{-1} dz.
\end{equation}
This projection allows us to isolate the invariant subspace associated with the spectral cluster $\Lambda_0$.

\begin{theorem}\label{theorem:Riesz_projection}
    The operator $P_\Gamma$ defined in \eqref{eq:Riesz_projection} is a~projection that commutes with $A$. Furthermore, its range $\mathcal{R}(P_\Gamma)$ is an invariant subspace of $A$, and the spectrum of $A$ restricted to $\mathcal{R}(P_\Gamma)$ is exactly $\Lambda_0$.
\end{theorem}

Following the methodology presented in \cite{Gavin2017, Polizzi2009}, the integral in \eqref{eq:Riesz_projection} can be evaluated numerically via quadrature. By computing an orthonormal basis $V \in \C^{n\times d}$ for $\mathcal{R}(P_\Gamma)$ and performing a~Rayleigh-Ritz projection, we obtain the reduced matrix $\tilde{A} := V^{H}AV \in \C^{d \times d}$. This procedure reduces the original problem to a~$d$-dimensional matrix whose spectrum consists exactly of the cluster $\Lambda_0$.

Shifting the spectrum of the matrix, so that its center is in the origin, we reduce the problem to matrices having only zero eigenvalue, i.e., nilpotent matrices. Recall that if $A$ is such matrix, then the number of Jordan blocks of size at least $k$ is given by: 
    \begin{equation}\label{eq:N}
        \mathrm{rank}\left(A^{k-1}\right) - \mathrm{rank}\left(A^k\right).
    \end{equation}
Consequently, the sequence of ranks of successive powers of $A$ uniquely determines the sizes of the Jordan blocks \cite{Golub1976}. Due to the nilpotency of $A$, in non-perturbed case always $A^d = 0$.

\section{Transformer-based Model for Recognition of Largest Attainable Jordan-block Size}

Code for reproducibility of conducted experiments is available online at \href{https://github.com/Michozord/Jordan_Transformer}{{https://github.com/Michozord/Jordan\_Transformer}}.

\subsection{Machine-learning framework}
Our goal is to estimate the largest possibly attainable Jordan block in some neighbourhood of a~given matrix. Following the method presented in Section \ref{section:nilpotent}, we can assume the matrix $A \in \R^{d \times d}$ to have exactly one cluster of eigenvalues. We restricted our model to real-valued matrices, as they are relevant for the applications in dynamical systems and engineering. 

Thus, we consider matrices $A \in \R^{d \times d}$ of the form 
\begin{equation}\label{eq:T_form} 
    a~= S\left(J+\eps E\right)S^{-1}
\end{equation}
where $J$ is a~Jordan matrix with one eigenvalue $\lambda = 0$, and $S$~is sampled from a~multivariate standard normal distribution subject to the constraint $\kappa(S) \leqslant 200d$. Using the approximation for the distribution of $\kappa(S)$ from Theorem \ref{theorem:Edelman}, we cover 99\% of matrices of $d\times d$ real Ginibre matrices with i.i.d. $\mathcal{N}(0,1)$ entries. We chose this boundary as it represents sufficiently wide class of similarity matrices to pass tests on matrices of different forms (see Section \ref{section:ood_tests}). On the other hand, our model struggles with similar problems as numerical approaches to related problems. Therefore, we could not allow Jordan matrices of arbitrary high conditioning of the Jordan basis.

The perturbation matrix is of the form $\eps E$, where $\eps>0$ is the perturbation magnitude and $E$ is a~normalized real-valued Ginibre matrix with i.i.d. entries sampled from $\mathcal{N}(0, 1/d)$. Thus, $\|\eps E\| = \eps$. Our objective is to predict the size $m$ of the largest block in $J$ given the matrix $A$.

\begin{figure}[tbp]
    \centering
    \includegraphics[width=\linewidth]{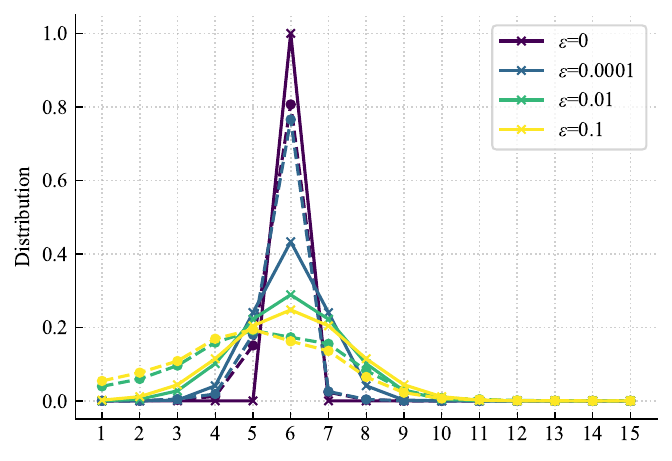}
    \caption{Comparison of the predicted distributions $\hat{q}$ (dashed lines) with the soft-target distributions $q_{m_{\mathrm{true}}, \eps}$ (solid lines) for matrices with $d=15$, $m_{\mathrm{true}}=6$, and varying $\eps$.} 
    \label{fig:JordanNet_q_example}
\end{figure}

For a~ground-truth category $m_{\mathrm{true}} \in \{1, \dots d\}$ and a~perturbation magnitude $\eps > 0$, we define a~soft-target distribution $q_{m_{\mathrm{true}}, \eps}$ over the classes $\{1, \dots d\}$. If the perturbation is negligible ($\eps \leqslant \eps_0$), the distribution is a~one-hot encoding; otherwise, it follows a~discrete Gaussian-like kernel:
\begin{equation}\label{eq:soft-target}
    q_{m_{\mathrm{true}}, \eps}(i) = \begin{cases}
        \delta_{i, m_{\mathrm{true}}}, & \text{if} \, \eps \leqslant \eps_0, \\
        \frac{\exp\left( - \frac{(i - m_{\mathrm{true}})^2}{2\tau^2} \right)}{\sum_{j=1}^d \exp\left( - \frac{(j - m_{\mathrm{true}})^2}{2\tau^2} \right)}, & \text{otherwise}
    \end{cases},
\end{equation}
where the temperature $\tau$ scales logarithmically with the perturbation:
\begin{equation}
\tau = C \log\left(1 + \frac{\eps}{\eps_0}\right).
\end{equation}
We set $\eps_0 = 10^{-8}$ and $C=0.1$. This formulation ensures that as $\eps$ increases, the target distribution diffuses symmetrically around $m_{\mathrm{true}}$, reflecting higher uncertainty for larger perturbations. Conversely, as $\eps \to \eps_0^+$, the distribution sharpens toward the true class. An example of the $q_{m_{\mathrm{true}}, \eps}$ distribution is presented in Figure \ref{fig:JordanNet_q_example}.

\subsection{Model}
In our initial experiments on low-dimensional matrices, we employed a~simple MLP. As input, we considered a~matrix $A$ containing a~single eigenvalue cluster, without any prior feature extraction. We trained the model using accuracy as the quality metric and cross-entropy as the loss function. The results convinced us that the problem can be learned by a~neural network-based model; nevertheless, this architecture was too simple to achieve satisfactory accuracy and could not be effectively scaled to higher dimensions. Consequently, we developed a~more sophisticated architecture. 

The workflow of the model presented in this paper is motivated by the properties of nilpotent matrices quantified in Section \ref{section:nilpotent}. As input, we take a~sequence $\left(A, A^2, \dots, A^{d}\right)$. Note that in a~non-perturbed case, $A^{d}$ is already zero matrix, regardless of the structure of $J$. However, we mainly focus on matrices with some perturbation, where matrix $A^d$ does not fully vanish and can indicate internal perturbation magnitude. Our experiments showed that considering one additional element of the sequence improves recognition of large blocks. The relevance of this feature is discussed further in Section \ref{section:attention}.

To ensure numerical stability of the matrix multiplication, we perform a~doubling test: we check whether $\|\widehat{A^{2k}} - (\widehat{A^k})^2\| \leqslant 10^{-12}$ for $k=1,\dots, \lfloor{d}/2\rfloor$, where $\widehat{A^k}$ denotes computed value of $A^k$. We discarded matrices violating this condition from the test set. Note that the condition $\rho(A) \leqslant 1$ is not only a~boundary for the size of the eigenvalue cluster we consider, but ensures numerical stability and avoids exploding eigenvalues of matrices $A^k$. 

We plug the input matrices into an encoder that embeds them into a~lower-dimensional latent space, process those embeddings with a~Transformer and finally use a~lightweight classifier to produce a~vector of logits representing possible maximal Jordan block sizes, $m\in \{1,\dots,d\}$.

While the task is mathematically straightforward for a~"perfect" matrix $A_0$ with spectral radius $\rho(A_0)=0$, we specifically consider cases with numerical perturbations, where the input $A^k$ represents a~perturbed version of the idealized $A_0^k$.

Encoders perform well at extracting robust features and input denoising tasks \cite{Vincent2010, Kulkarni2023}. Therefore, we employ an Encoder $E_d : \R^{d \times d} \rightarrow \R^{d_\mathrm{latent}}$ for the processing of the input sequence. We embed the sequence into a~latent space $\R^{d_\mathrm{latent}}$ ($d_\mathrm{latent}=32$) using a~multilayer perceptron consisting of three hidden layers with $d_{\mathrm{hidden}}=128$ neurons each and ReLU activations. We denote the resulting embedding as $y_k := E_d(A^k)$. While all matrices in a~sequence share the same encoder, each dimension $d$ requires a~dedicated encoder. Encoder outputs are layer-normalized \cite{ba2016}.

At the core of the model, we employ a~Transformer architecture \cite{Attention}, which is well-suited for processing sequences of vectors. Our implementation uses 2 layers, feed-forward networks of dimension $d_{\mathrm{hidden}}$, and 4 attention heads. Following the Transformer layers, we apply mean pooling to obtain $\bar{z} = 1/d \sum_{i=1}^d z_i \in \R^{d_{\text{latent}}}$, where $(z_1, \dots, z_d)$ is the Transformer output. Notably, the Transformer is independent of the matrix dimension $d$ and is shared across all considered dimensions. 

\begin{figure}[tbp]
    \centering
    \includegraphics[width=3in]{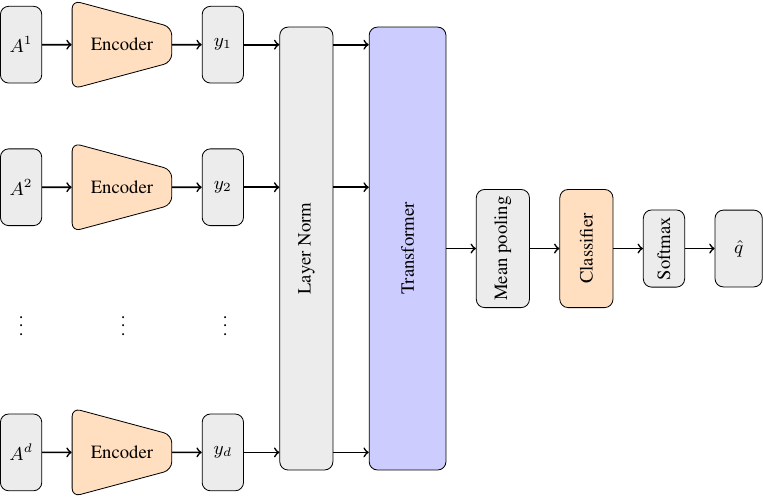}
    \caption{Model architecture. Orange modules are dedicated to specific dimension, while the blue Transformer is dimension-independent.}
    \label{fig:JordanNet}
\end{figure}

The aggregate vector $\bar{z}$ is finally passed to a~dimension-specific Classification Head $C_d : \R^{d_{\text{latent}}} \rightarrow \R^d$, which is a~multilayer perceptron with a~single 128-dimensional hidden layer. The output is a~$d$-dimensional vector of logits; applying the softmax function yields the predicted distribution $\hat{q}$ over all possible block sizes. The model architecture is shown in Figure \ref{fig:JordanNet}, and Table \ref{table:JordanNet_weights} summarizes the parameter counts.

\begin{table*}[tbp]
    \caption{Parameter counts for model components.}
    \label{table:JordanNet_weights}
    \centering
    \setlength{\tabcolsep}{4pt}
    \renewcommand{\arraystretch}{1.3}
    \begin{tabular}{llr}
        \toprule
        \textbf{Module} & \textbf{General Formula ($d, \dhidden, \dlatent$)} & \textbf{Weights ($\dhidden=128, \dlatent=32$)} \\ 
        \midrule
        Encoder     & $\dhidden(d^2 + 2\dhidden + \dlatent + 3) + \dlatent$ & $128d^2 + 37,280$ \\ 
        Transformer & $24\dlatent^2 + 26\dlatent$                           & $25,408$          \\
        Classifier  & $\dhidden(\dlatent + d + 1) + d$                      & $129d + 4,224$    \\ 
        \midrule
        \textbf{Total} & \multicolumn{2}{r}{%
            $\underbrace{25,408}_{\text{Shared}} + \underbrace{128d^2 + 129d + 41,504}_{\text{Per dimension } d}$%
        } \\ 
        \bottomrule
    \end{tabular}
\end{table*}

\subsection{Training}
Training was split into two phases. In the first phase, we trained the Encoders $E_d$, Heads $C_d$, and the Transformer simultaneously for dimensions $d = 4, 6, 9, 12, 15, 28$. This multi-dimension training is essential for ensuring the Transformer generalizes across different matrix sizes. In the second phase, we trained additional Encoders and Classifiers for additional dimensions $d$ (also exceeding dimensions involved in Transformer training) while freezing the Transformer parameters. This was done to test the Transformer's universality and the cost of extending the model to new dimensions, including those exceeding the maximum dimension of the first phase.

In Phase 1, the dataset consisted of matrices generated according to \eqref{eq:T_form}. The perturbation magnitude $\eps$ was sampled independently for each matrix using a~spike-and-log-uniform distribution: $\eps = 0$ with probability $0.1$, and $\log\eps \sim \mathcal{U}(\log \eps_{\min}, \log \eps_{\max})$ otherwise, with $\eps_{\min} = 10^{-16}$ and $\eps_{\max} = 0.1$. To ensure fulfilling the assumption $\rho(A) \leqslant 1$, if the sampled matrix $A$ violated it, we normalized $A$, dividing it by $\rho(A)$ and appropriately scaling $\eps$ in the soft-target.

This dataset included 148,000 matrices (2,000 per class $m$ per dimension). The Phase 2 datasets contained 1,200 matrices per class per dimension, generated in the same manner.

For both phases, 20\% of the data was reserved for validation. The loss function was the average Kullback-Leibler (KL) divergence between the soft-target $q_{m_{\mathrm{true}}, \eps}$ and the predicted distribution $\hat{q}$: 
\begin{equation}
    D_{KL}(q_{m_{\mathrm{true}}, \eps} \parallel \hat{q}) := \sum_{m=1}^d q_{m_{\mathrm{true}}, \eps}(m) \log \frac{q_{m_{\mathrm{true}}, \eps}(m)}{\hat{q}(m)}.
\end{equation} 
We used the Adam optimizer \cite{Adam2017} with cosine annealing \cite{Loshchilov2017}, a~maximum of 40 epochs, no restarts, and a~learning rate range of $\eta \in [0, 5\cdot10^{-4}]$. Training stopped if the validation loss failed to improve by more than $10^{-4}$ for 8 consecutive epochs. Phase 1 concluded at 27 epochs.

\section{Results}

For the model's output $\hat{q}$, let $\hat{m} := \mathrm{argmax}_m\hat{q}(m)$ denote the predicted Jordan block size. We evaluated the model using three quality metrics: average KL divergence, which is consistent with the training loss; model accuracy, defined as the proportion of matrices where $m_{\mathrm{true}} = \hat{m}$; and relaxed accuracy, defined as the proportion where $|m_{\mathrm{true}} - \hat{m}| \leqslant 1$. 

\subsection{In-distribution tests}
We generated test sets containing matrices from the same distribution as the training sets, i.e., of the form \eqref{eq:T_form} with varying perturbation magnitudes $\eps \in [0, 0.1]$. For each dimension, we generated 10,000 matrices per class in the same manner as the training set. Tables \ref{tab:metrics_eps} and \ref{tab:metrics_rho} present the obtained accuracy scores depending on the perturbation magnitude $\eps$ and the spectral radius $\rho(A)$, respectively. 

Comparing metrics across dimensions from the first and second phases of training, we observe that the model achieves comparably high accuracy for dimensions not included in the Transformer training. This includes dimensions $d=33$ and $35$, which exceed the maximum dimension involved in the first training phase. This demonstrates the Transformer's adaptability to new dimensions and suggests that it has learned a~dimension-agnostic mathematical logic governing the Jordan structure, rather than simply memorizing patterns for specific matrix sizes. 

As expected, the model's performance drops as $\rho(A)$ approaches 1 or $\eps$ approaches $\eps_{\max}=0.1$. In these cases, the perturbations become larger, hindering the model from rendering correct predictions. Notably, the accuracy@$\pm1$ remains above $0.9$ even for $\rho(A)$ up to 0.75 and $\eps$ reaching $\eps_{\max} = 0.1$. This indicates that the model typically misclassifies the matrix by at most 1. This case is exemplified in Figure \ref{fig:JordanNet_q_example}, where for $\eps=0.01$ and $\eps=0.1$, the model's predicted distribution $\hat{q}$ is shifted by one to the left of the true value.

\begin{table*}[t]
\centering
\caption{Performance metrics across dimensions $d$ and perturbation magnitude $\eps$. $\dagger$ denotes dimensions involved in training of the Transformer.}
\label{tab:metrics_eps}

\setlength{\tabcolsep}{4pt}
\begin{tabular}{llccccccccc}
\toprule
 & $\eps$
 & $d=6^\dagger$ & $d=9^\dagger$ & $d=12^\dagger$ & $d=15^\dagger$ & $d=19$
 & $d=25$ & $d=28^\dagger$ & $d=33$ & $d=35$ \\
\midrule

\multirow{4}{*}{\textbf{Avg.\ KL div.}}
& $0$
    & 0.042 & 0.035 & 0.010 & 0.006 & 0.032
    & 0.019 & 0.002 & 0.011 & 0.024 \\
& $(0,10^{-3}]$
    & 0.149 & 0.119 & 0.104 & 0.098 & 0.188
    & 0.200 & 0.082 & 0.108 & 0.134 \\
& $(10^{-3},10^{-2}]$
    & 0.010 & 0.005 & 0.009 & 0.005 & 0.032
    & 0.082 & 0.007 & 0.047 & 0.050 \\
& $(10^{-2},10^{-1}]$
    & 0.030 & 0.032 & 0.031 & 0.029 & 0.095
    & 0.175 & 0.046 & 0.125 & 0.143 \\

\midrule

\multirow{4}{*}{\textbf{Acc}}
& $0$
    & 0.999 & 1.000 & 1.000 & 1.000 & 1.000
    & 1.000 & 1.000 & 1.000 & 1.000 \\
& $(0,10^{-3}]$
    & 0.999 & 1.000 & 1.000 & 1.000 & 0.990
    & 0.993 & 1.000 & 0.996 & 0.993 \\
& $(10^{-3},10^{-2}]$
    & 0.986 & 0.990 & 0.987 & 0.993 & 0.905
    & 0.886 & 0.992 & 0.946 & 0.947 \\
& $(10^{-2},10^{-1}]$
    & 0.841 & 0.887 & 0.879 & 0.897 & 0.750
    & 0.753 & 0.892 & 0.767 & 0.784 \\

\midrule

\multirow{2}{*}{\textbf{Acc@$\pm1$}}
& $(10^{-3},10^{-2}]$
    & 1.000 & 1.000 & 0.999 & 1.000 & 0.999
    & 0.996 & 1.000 & 0.997 & 0.997 \\
& $(10^{-2},10^{-1}]$
    & 0.980 & 0.982 & 0.985 & 0.985 & 0.955
    & 0.938 & 0.977 & 0.941 & 0.943 \\

\bottomrule
\end{tabular}
\end{table*}
\begin{table*}[t]
\centering
\caption{Performance metrics across dimensions $d$ and spectral radius $\rho(A)$. $\dagger$ denotes dimensions involved in training of the Transformer.}
\label{tab:metrics_rho}

\setlength{\tabcolsep}{4pt}
\begin{tabular}{llccccccccc}
\toprule
 & $\rho(A)$
 & $d=6^\dagger$ & $d=9^\dagger$ & $d=12^\dagger$ & $d=15^\dagger$ & $d=19$
 & $d=25$ & $d=28^\dagger$ & $d=33$ & $d=35$ \\
\midrule

\multirow{5}{*}{\textbf{Avg.\ KL div.}}
& $[0,10^{-8}]$
    & 0.042 & 0.035 & 0.010 & 0.007 & 0.033
    & 0.021 & 0.003 & 0.012 & 0.025 \\
& $(10^{-8},0.25]$
    & 0.120 & 0.103 & 0.086 & 0.080 & 0.101
    & 0.151 & 0.069 & 0.138 & 0.124 \\
& $(0.25,0.5]$
    & 0.030 & 0.032 & 0.067 & 0.080 & 0.218
    & 0.138 & 0.052 & 0.068 & 0.057 \\
& $(0.5,0.75]$
    & 0.022 & 0.043 & 0.035 & 0.025 & 0.115
    & 0.232 & 0.068 & 0.080 & 0.143 \\
& $(0.75,1]$
    & -- & 0.081 & 0.053 & 0.082 & 0.211
    & 0.292 & 0.076 & 0.170 & 0.183 \\

\midrule

\multirow{5}{*}{\textbf{Acc}}
& $[0,10^{-8}]$
    & 0.999 & 1.000 & 1.000 & 1.000 & 1.000
    & 1.000 & 1.000 & 1.000 & 1.000 \\
& $(10^{-8},0.25]$
    & 0.984 & 0.995 & 0.994 & 0.996 & 0.992
    & 0.996 & 0.999 & 0.987 & 0.993 \\
& $(0.25,0.5]$
    & 0.879 & 0.962 & 0.984 & 0.994 & 0.974
    & 0.994 & 0.999 & 0.981 & 0.995 \\
& $(0.5,0.75]$
    & 0.837 & 0.805 & 0.857 & 0.934 & 0.805
    & 0.895 & 0.994 & 0.979 & 0.972 \\
& $(0.75,1]$
    & -- & 0.591 & 0.489 & 0.596 & 0.553
    & 0.512 & 0.810 & 0.721 & 0.737 \\

\midrule

\multirow{2}{*}{\textbf{Acc@$\pm1$}}
& $(0.5,0.75]$
    & 0.998 & 0.971 & 0.981 & 0.989 & 0.986
    & 0.994 & 0.999 & 0.998 & 0.999 \\
& $(0.75,1]$
    & -- & 0.955 & 0.983 & 0.945 & 0.864
    & 0.863 & 0.959 & 0.921 & 0.932 \\

\bottomrule
\end{tabular}
\end{table*}

To evaluate the versatility and scalability of the architecture, we tested the robustness of a~model trained on matrices of dimension $d_1$ when applied to matrices of a~smaller dimension $d_2 < d_1$ by padding the unused entries with zeros. Specifically, we generated matrices $A_0 \in \R^{d_2 \times d_2}$ using the method described above and fed the model matrices of the form:
\begin{equation}\label{eq:oplus_matrices}
    a~= A_0 \oplus \textbf{0}_{(d_1 - d_2)\times(d_1-d_2)} \in \R^{d_1 \times d_1}.
\end{equation}
Note that this transformation does not impact the true class $m_{\mathrm{true}}$ or the distance $\eps$ to the defective matrix with the largest attainable block. These results are presented in Table \ref{tab:oplus_results}. 

\begin{table}[t]
\centering
\caption{Performance metrics for matrices of the form \eqref{eq:oplus_matrices} across perturbation magnitude $\eps $}
\label{tab:oplus_results}

\setlength{\tabcolsep}{4pt}
\renewcommand{\arraystretch}{1.2}
\begin{tabular}{l l c c c}
\toprule
 & $\eps$
 & \shortstack{$d_2=16$\\$d_1=19$} 
 & \shortstack{$d_2=27$\\$d_1=28$} 
 & \shortstack{$d_2=29$\\$d_1=33$} \\
\midrule

\multirow{4}{*}{{\textbf{Avg.\ KL div.}}}
& $0$
    & 0.029 & 0.003 & 0.014 \\
& $(10^{-12},10^{-3}]$
    & 0.260 & 0.104 & 0.327 \\
& $(10^{-3},10^{-2}]$
    & 0.127 & 0.007 & 0.425 \\
& $(10^{-2},10^{-1}]$
    & 1.114 & 0.140 & 2.345 \\

\midrule

\multirow{4}{*}{{\textbf{Acc}}}
& $0$
    & 1.000 & 1.000 & 1.000 \\
& $(10^{-12},10^{-3}]$
    & 0.982 & 1.000 & 0.999 \\
& $(10^{-3},10^{-2}]$
    & 0.801 & 0.994 & 0.941 \\
& $(10^{-2},10^{-1}]$
    & 0.507 & 0.855 & 0.549 \\

\midrule

\multirow{2}{*}{{\textbf{Acc@$\pm1$}}}
& $(10^{-3},10^{-2}]$
    & 0.978 & 0.999 & 0.971 \\
& $(10^{-2},10^{-1}]$
    & 0.764 & 0.944 & 0.690 \\

\bottomrule
\end{tabular}
\end{table}

\begin{figure}[!tbh]
\centering
\subfloat[]{\includegraphics[width=\linewidth]{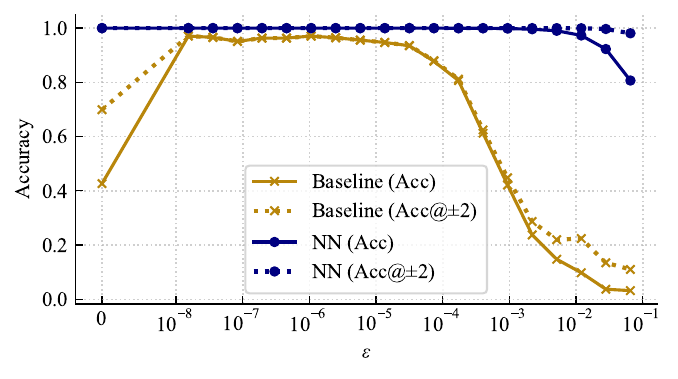}%
\label{fig_first_case}}
\hfil
\subfloat[]{\includegraphics[width=\linewidth]{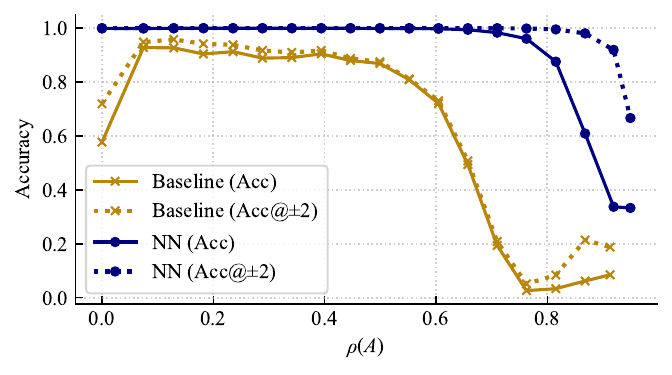}%
\label{fig_second_case}}
\caption{Comparison of accuracy and accuracy@$\pm2$ for the naive numerical classifier (Baseline) and our model (NN) for $d=28$ as a~function of (a) perturbation magnitude $\eps$ and (b) spectral radius $\rho(A)$.}
\label{fig:naive_baseline}
\end{figure}

To provide stronger evidence of our model's robustness against perturbations, we compare its performance to a~baseline numerical algorithm based on the identity \eqref{eq:N}. For a~matrix $A$ of the form \eqref{eq:T_form}, this baseline algorithm computes the rank of consecutive powers of $A$ via the SVD-based algorithm presented in \cite{Press2007} (the standard implementation for matrix rank in Matlab and NumPy). Figure \ref{fig:naive_baseline} presents the results. Our model consistently outperforms the baseline algorithm, especially for larger perturbation magnitudes. This demonstrates that the proposed Encoder-Transformer architecture effectively mitigates heavy noise, leading to correct predictions even in highly perturbed cases.
 
\subsection{Out-of-distribution tests}\label{section:ood_tests}
To evaluate the performance of our model on matrices from a~different distribution than the one used during training and to ensure its versatility, we constructed test matrices with selected Segre characteristics (i.e. sequence of sizes of Jordan blocks corresponding to the zero eigenvalue) for dimensions $d = 8, 12, 13$, where the matrices of size 8 were cast to $d=9$ using \eqref{eq:oplus_matrices}. We tested matrices of the form $A = A_0 + \eps E$, where $E$ is a~real Ginibre matrix with i.i.d. entries scaled to $\|E\|=1$, and $A_0$ is a~matrix with a~single eigenvalue at the origin possessing the desired Segre characteristics. This method of constructing test matrices originates from \cite{Ran2010}. The $\eps$ varied between $10^{-7}$ and $10^{-3}$, and the corresponding $\rho(A)$ is reported in the table. Note that, unlike the notation in the remainder of this article, here $\eps=\delta_k(A)$ refers to the distance to the defective matrix. 

The results for selected characteristics are presented in Table \ref{tab:segre_metrics}. This experiment demonstrates that our model also performs well on matrices constructed differently from those in the training set.

\begin{table}[t]
\centering
\caption{Performance on out-of-distribution matrices with different Segre characteristics and spectral radii}
\label{tab:segre_metrics}
\setlength{\tabcolsep}{4pt}

\begin{tabular}{lcccc}
\toprule
 & (3,2,2,1) & (2,2,2,2) & (3,3,3,3,1) & (4,3,3,2) \\
\midrule

\multicolumn{5}{l}{\textbf{Acc}} \\
$\rho(A) \in [0, 0.01)$
    & 0.658 & 0.999 & 0.995 & 0.247 \\
$\rho(A) \in [0.01, 0.03)$
    & 0.806 & 1.000 & 0.997 & 0.687 \\
$\rho(A) \in [0.03, 1)$
    & 0.967 & 1.000 & 0.949 & 0.804 \\
\midrule

\multicolumn{5}{l}{\textbf{Acc@$\pm1$}} \\
$\rho(A) \in [0, 0.01)$
    & 0.990 & 1.000 & 1.000 & 0.929 \\
$\rho(A) \in [0.01, 0.03)$
    & 0.997 & 1.000 & 1.000 & 0.985 \\
$\rho(A) \in [0.03, 1)$
    & 1.000 & 1.000 & 1.000 & 0.996 \\
\bottomrule
\end{tabular}
\end{table}

\section{Inspecting The Model}\label{section:attention}
\begin{figure}[tbp]
    \centering
    \includegraphics[width=\linewidth]{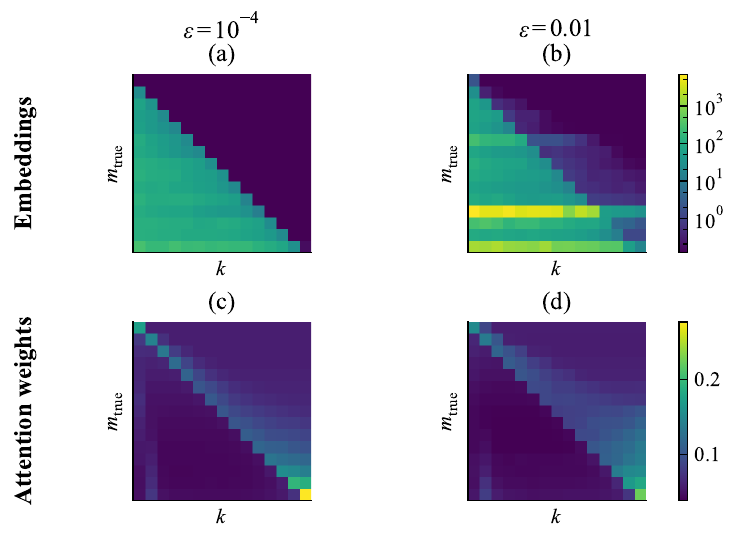}
    \caption{Average internal data of the model for $d=15$ across a~dataset of 100 matrices per maximal block size with perturbation magnitudes $\eps=10^{-3}$ and $\eps=0.1$. (a)--(b): Norm of embedding $\|y_k\|$ vs. ground-truth category $m_{\mathrm{true}}$ (asinh scale). (c)--(d): Importance of embeddings $y_k$ defined as the average of all attention weights corresponding to the element of the Transformer input vs. ground-truth $m_{\mathrm{true}}$.}
    \label{fig:Attention}
\end{figure}

To analyse the model's internal workflow, we examine the learned internal latent representations and attention weights. As expected, the embeddings $y_k\in \R^{\dlatent}$, $k=1,\dots,d-1$, produced by the encoder $E_d$, reflect the nullity of matrices $A^k$ through their magnitude, i.e., $\|y_k\| \approx 0$ if and only if $A^k = 0$ in the noiseless, theoretical case. This relation is depicted in Figure \ref{fig:Attention}(a)--(b).

Even though it may appear to be a~straight-forward task, we experimented with shallower networks for this role but observed an immediate degradation in recognition accuracy. This is evidence of the higher complexity of algebraic features required for Jordan structure recognition compared to the ViT task, where a~linear projection of an image patch is sufficient \cite{Kolesnikov2021}. The encoders effectively act as implicit denoisers, extracting algebraic features from the input sequence and mapping them into a~low-dimensional latent space that is universal across dimensions. At a~higher perturbation magnitude $\eps=0.1$, we observe the influence of noise even in cases where theoretically $A^k=0$. This is visible in the upper triangle of the plot in Figure \ref{fig:Attention}(b).

The Transformer performs dimension-agnostic reasoning over the embeddings of nilpotent matrices. In Figure \ref{fig:Attention}(c)--(d), we present the importance of each input $y_k$ as the mean of the attention weights attributed to each token. The model tends to put the highest attention on the token corresponding to the first vanishing $A^k$ matrix and its successors. In context of \eqref{eq:N}, this behaviour is intuitive. These tokens contain relevant information where the transition to near-zero matrices occurs. Here, the effect of the noise in the input matrix $A$ is most visible, especially by higher perturbations. For the same reason, we plug in matrix $A^d$ to the model, which allows the model to recognize the perturbation magnitude.

\section{Conclusion}
In this paper, we introduced a~Transformer-based model for determining the maximal attainable Jordan block size under a~boundary constraint on the conditioning of the Jordan basis, $\kappa(S) \leqslant 200d$. Our architecture, although inspired by Vision Transformers designed for image recognition, successfully overcomes the complete lack of spatial dependence in the target matrices. This was achieved by using algebraically motivated features, namely the gradually vanishing powers of the objective matrix. Introducing the noise-aware soft-target training distribution allowed the model to gracefully handle numerical uncertainties caused by noise and the underlying specific nature of the Jordan structure itself. Thanks to this approach, the model misclassifies the matrix by at most $\pm1$ in most cases, which we have observed across multiple conducted experiments.

The successful training of additional encoders and classifiers with frozen Transformer weights, alongside the inspection of the attention weights, shows that the model robustly learned dimension-agnostic algebraic rules. While numerical approaches rely on similar algebraic properties of nilpotent matrices, our model exhibited significantly higher robustness against numerical noise. 

The presented model still exhibits some limitations: we considered only real-valued matrices and restricted the problem to matrices having a~single eigenvalue cluster, e.g., after applying a~Riesz projection. Finally, the generalization of the presented methods to matrix pencils remains an open problem for future research.

\bibliographystyle{IEEEtran}
\bibliography{bib}

\end{document}